\documentclass[11pt,letterpaper,reqno]{amsart}

\usepackage[margin=1in]{geometry}
\usepackage{amsmath,amssymb,amsthm,mathtools}
\usepackage{hyperref}
\usepackage[nameinlink,capitalize]{cleveref}
\usepackage{enumitem, comment, amssymb}
\usepackage{color,soul}

\newtheorem{theorem}{Theorem}[section]

\newtheorem{lemma}[theorem]{Lemma}
\newtheorem{corollary}[theorem]{Corollary}
\newtheorem{example}[theorem]{Example}

\usepackage{scalerel}[2014/03/10]
\usepackage[usestackEOL]{stackengine}
\def\intavg{\,\ThisStyle{\ensurestackMath{%
    \stackinset{c}{0\LMpt}{c}{0\LMpt}{\SavedStyle-}{\SavedStyle\phantom{\int}}}%
    \setbox0=\hbox{$\SavedStyle\int\,$}\kern-\wd0}\int}

\theoremstyle{definition}
\newtheorem{definition}[theorem]{Definition}
\newtheorem{remark}[theorem]{Remark}
\newtheorem{question}[theorem]{Question}

\title{Regularity of Weak Solutions near Singular Minimizing Cones for the One-Phase Free Boundary Problem}

\author{Minyoung Jeong}
\address{Department of Mathematics, University of California San Diego}
\email{m1jeong@ucsd.edu}


\begin{document}

\begin{abstract}
Alt and Caffarelli \cite{AC81} proved an Allard-type $\epsilon$-regularity theorem near the flat one-homogeneous minimizer of the one-phase free boundary problem. A natural follow-up question is whether an analogous regularity result holds near singular minimizing cones. In this paper, we prove such an $\epsilon$-regularity theorem near a strongly integrable one-phase minimizing cone with an isolated singularity. More precisely, we show that, after a suitable translation and rotation, the free boundary of any non-degenerate weak solution sufficiently close to the cone is graphical over the free boundary of a leaf of the De Silva--Jerison--Shahgholian foliation \cite{DSJS22}. This result is analogous to recent $\epsilon$-regularity theorems in minimal surface theory by Edelen \cite{E24} and by Edelen--Spolaor \cite{ES23}; our proof adapts the strategy of Edelen \cite{E24} to the one-phase setting.
\end{abstract}

\maketitle
\markboth{}{}
\tableofcontents

\section{Introduction}\label{Introduction}
The one-phase free boundary problem is a prototypical model for stationary free interfaces arising in phase-transition and fluid-interface problems \cite{VW11, WZ12}. In this model, the interface is not prescribed a priori, but is instead determined as part of the variational problem. Mathematically, given an open set $\Omega \subseteq \mathbb R^d$, we consider the Alt-Caffarelli energy
\[
\mathcal{F}(u; \Omega) := \int_{\Omega} |\nabla u|^2 + \chi_{\{u > 0\}} \, dx
\]
among $u \in H^1(\Omega)$ satisfying $u \geq 0$ in $\Omega$ and $u = g$ on $\partial \Omega$.

Formally, a minimizer $u$ of this energy solves the following Euler-Lagrange equation
\begin{equation}
  \label{eq:serrin-prob}
  \left\{
  \begin{aligned}
    \Delta u &= 0 \quad&&\text{in } \{u > 0\} \cap \Omega\\
    u &= 0 \quad&&\text{on } \partial \{u > 0\} \cap \Omega \\
    \partial_{\nu}u &= 1 \quad&&\text{on } \partial \{u > 0\} \cap \Omega.
  \end{aligned}
  \right.
\end{equation}
where $\nu$ denotes the inner unit normal to the free boundary $\partial \{u > 0\} \cap \Omega$.

Alt and Caffarelli \cite{AC81} proved that minimizers of $\mathcal{F}$ are locally Lipschitz in $\Omega$ and that the free boundary $\partial\{ u > 0\} \cap \Omega$ can be decomposed into regular and singular parts, where the regular set is a $C^{1, \alpha}$ hypersurface for some $\alpha > 0$. One of their main tools is an Allard-type $\epsilon$-regularity theorem near the flat one-homogeneous minimizer: if a blow-up of $u$ at a free boundary point $x_0$ is, up to rotation, equal to $(x_d)^+$, then the free boundary is a $C^{1,\alpha}$ hypersurface in a neighborhood of $x_0$.

Weiss \cite{W98, W99} subsequently proved a monotonicity formula which implies that blow-up limits at free boundary points are one-homogeneous. This leads naturally to the following question:
\begin{question} \label{question01}
Suppose that a blow-up of a minimizer $u$ at a free boundary point $x_0$ is a singular minimizing cone $U$. What can be said about the structure of the free boundary near $x_0$?
\end{question}
More generally, one may ask the following question:
\begin{question} \label{question02}
Suppose that $u$ is a minimizer, or more generally a weak solution of the one-phase problem, which is sufficiently close to a singular minimizing cone $U$. What can be said about the structure of the free boundary of $u$ near $x_0$?
\end{question}

There has been substantial progress in understanding singularities of minimizing free boundaries. Caffarelli, Jerison, and Kenig \cite{CJK04} proved that minimizers have no singular free boundary points in dimensions $d \leq 3$, and Jerison and Savin \cite{JS15} extended this result to $d \leq 4$. On the other hand, De Silva and Jerison \cite{DSJ09} constructed an explicit one-homogeneous global minimizer with an isolated singularity in dimension $d = 7$. These results established that the critical dimension $d^*$, defined as the smallest dimension admitting a singular minimizing cone, satisfies $5 \leq d^* \leq 7$. A recent preprint by OpenAI \cite{OpenAI26} reports a proof that $d^* = 7$.

Moreover, considerable work has been devoted to the size and geometric structure of the singular set. Using quantitative stratification and the rectifiable-Reifenberg theory, Edelen and Engelstein \cite[Theorem 1.1]{EE19} proved that the singular set of a minimizer is countably $(d-d^*)$-rectifiable and has locally finite $\mathcal{H}^{d-d^*}$-measure. More recently, Kriventsov and Weiss \cite{KW25} studied a broader class of non-minimizing variational solutions. They proved that the entire topological free boundary, including degenerate singular points, is countably $(d-1)$-rectifiable and has locally finite $\mathcal{H}^{d-1}$-measure.

Several recent works are particularly relevant to the present paper. Given a one-homogeneous global minimizer $U$ with an isolated singularity at the origin, De Silva, Jerison, and Shahgholian \cite[Theorem 1.1]{DSJS22} constructed a smooth foliation around $U$:
\begin{theorem}[\cite{DSJS22}]
Let $U : \mathbb R^d \rightarrow [0, \infty)$ be a one-homogeneous global minimizer so that $\partial \{U > 0\} \setminus \{0\}$ is smooth. Then there exist global minimizers $\overline{U}, \underline{U} : \mathbb R^d \rightarrow [0, \infty)$
\begin{enumerate}
    \item $\partial \{\overline{U} > 0 \}$ and $\partial \{\underline{U} > 0 \}$ are analytic hypersurfaces in $\mathbb R^d$.
    
    \item $\underline{U} \leq U \leq \overline{U}$ and $\mathrm{dist}(\partial \{\overline{U} > 0 \}, \{0\}) = \mathrm{dist}(\partial \{\underline{U} > 0 \}, \{0\}) = 1$
    
    \item For $t > 0$, define
    \[
        \overline{U}_t := \frac{\overline{U}(t\cdot)}{t}, \qquad \underline{U}_t := \frac{\underline{U}(t\cdot)}{t}.
    \]
    Then the (positive) graphs of $(\overline{U}_t)_{t >0}, (\underline{U}_t)_{t >0}$, and $U$ foliate the half-space $\mathbb R^d \times (0, \infty)$, i.e. if we let $\overline{G}_t : = \{(x, \overline{U}_t(x)) : \overline{U}_t(x) > 0 \}$, $\underline{G}_t : = \{(x, \underline{U}_t(x)) : \underline{U}_t(x) > 0 \}$, and $G := \{(x, U(x)) : U(x) > 0 \}$ for $t > 0$, then $\overline{G}_t, \underline{G}_t, G$ are pairwise disjoint and
    \[
        \mathbb R^d \times (0, \infty) = G \cup \bigcup_{t > 0} \left( \overline{G}_t \cup \underline{G}_t \right)
    \]
    
    \item Let $V$ be a global minimizer. If $V \geq U$, then either $V \equiv U$ or $V \equiv \overline{U}_t$ for some $t > 0$. Likewise, if $V \leq U$, then either $V \equiv U$ or $V \equiv \underline{U}_t$ for some $t > 0$.
\end{enumerate}
\end{theorem}

Engelstein, Spolaor, and Velichkov \cite{ESV20} answered Question \ref{question01} when the blow-up cone $U$ has an isolated singularity. They proved that if a blow-up $U$ of a minimizer $u$ at a free boundary point has an isolated singularity, then $U$ is the unique blow-up limit at $x_0$, with a logarithmic rate of convergence, and the free boundary of $u$ near $x_0$ is a $C^{1, \log}$ graph over the free boundary of $U$.

More recently, Engelstein, Restrepo, and Zhao \cite{ERZ25} proved an analogous uniqueness-of-blow-up theorem, with a Hölder rate of convergence, for non-degenerate weak solutions whose blow-up cone is an integrable cone with an isolated singularity. They also established a rigidity theorem for global weak solutions asymptotic at infinity to strongly integrable minimizing cones:
\begin{theorem}[\cite{ERZ25}] \label{Theorem_blowdown_rigidity}
Let $u : \mathbb R^d \rightarrow [0, \infty)$ be a $(c_1, c_2)$-non-degenerate weak solution in the sense of Definition \ref{def_non_degenerate_weak_sol}. Let $U : \mathbb R^d \rightarrow [0, \infty)$ be a strongly integrable minimizing cone (see \cite[Definition 2.5]{ERZ25}) such that $\partial \{U > 0\} \setminus \{0\}$ is smooth. Suppose that there exist $x_0 \in \mathbb R^d$ and $R_j \uparrow \infty$ such that
\[
    \frac{u(x_0 + R_j\cdot)}{R_j} \rightarrow U
\]
locally uniformly in $\mathbb R^d$.\\
Then, up to a translation, $u$ coincides with $U, \overline{U}_t$, or $\underline{U}_t$ for some $t > 0$.
\end{theorem}

The two uniqueness results have different scopes: \cite{ERZ25} applies to non-degenerate weak solutions, and the limiting cone is integrable but need not be minimizing. On the other hand, \cite{ESV20} applies to minimizing solutions, and the limiting cone is minimizing but no integrability assumption is required, and the convergence rate is logarithmic in general. More recently,  Carducci and Velichkov \cite{CV26} proved uniqueness of blow-ups for Lipschitz stationary solutions at points where one blow-up has an isolated singularity, without assuming integrability of the limiting cone.

\subsection{Main result}
Let $S$ be a smooth hypersurface in $\mathbb R^d$ with a smooth unit normal vector field $\nu_S : S \rightarrow \mathbb S^{d-1}$. For $\xi : S \rightarrow \mathbb R$, define the graph of $\xi$ over $S$ as
\[
    G_{S}(\xi) := \{x + \xi(x)\nu_S(x) : x \in S\}.
\]
Also, for $(a, t, q) \in \mathbb R^d \times \mathbb R \times SO(d)$, define
\begin{equation} \label{Definition_F}
    F_{a,t,q} := 
    \begin{cases}
        a + q(\partial \{ \overline{U}_\frac{1}{t} > 0\}) & t > 0\\
        a + q(\partial \{ U > 0\}) & t = 0\\
        a + q(\partial \{ \underline{U}_{\frac{1}{|t|}} > 0\}) & t < 0.
    \end{cases}
\end{equation}
We define the regular part of the free boundary as $F^*_{a, t, q} := F_{a, t, q} \setminus \{a\}$. Note that $F^*_{a, t, q} = F_{a, t, q}$ when $t \ne 0$. In this paper, we answer Question \ref{question02} for strongly integrable minimizing cones with isolated singularities:
\begin{theorem} \label{main_theorem}
Let $U$ be a strongly integrable minimizing cone with $\mathrm{Sing}(U) = \{0\}$, and let $c_1, c_2, c_3 > 0$, with $c_1 \leq c_2$. For every $\epsilon > 0$, there exists $\delta(\epsilon, c_1, c_2, c_3, U) > 0$ such that the following holds:\\
Let $u : B_3 \rightarrow \mathbb R$ be a $(c_1, c_2, c_3)$-non-degenerate weak solution satisfying:
\[
\|u - U\|_{L^2(B_2)} < \delta.
\]
Then there exists $(a, t, q) \in \mathbb R^d \times \mathbb R \times SO(d)$ and $\xi \in C^2(F_{a, t, q}^* \cap B_{\frac{1}{2}}(a))$ such that
\[
    |a| + |t| + |\mathrm{Id} - q| \leq \epsilon
\]
and
\begin{equation}
\begin{cases}
    \partial \{u > 0\} \cap B_{\frac{1}{4}}(a) = \overline{G_{F_{a, t, q}^* \cap B_{\frac{1}{2}}(a)}(\xi)} \cap B_{\frac{1}{4}}(a) \\

    |x-a|^{-1}|\xi| + |\nabla \xi| + |x-a||\nabla^2 \xi| \leq \epsilon \text{ for $x \in F_{a, t, q}^* \cap B_{\frac{1}{2}}(a)$}.
\end{cases}
\end{equation}

\end{theorem}

The De Silva-Jerison cone \cite{DSJ09} is a minimizing cone with an isolated singularity and is strongly integrable \cite[Theorem 2.9]{ERZ25}. Thus, Theorem \ref{main_theorem} applies to the De Silva-Jerison cone. More recently, Firester, Tsiamis, and Wang \cite[Theorem 1.4]{FTW26b} proved minimality for a subfamily of cones with bi-orthogonality symmetry, although strong integrability of the new minimizing examples remains to be verified. Carducci and Velichkov \cite{CV26b} constructed a minimizing cone in dimension 26 that is neither axially symmetric nor of Lawson type. Analogous results in minimal surface theory were proved by Edelen, Spolaor \cite{ES23} and Edelen \cite{E24}. Our proof adapts the strategy of Edelen \cite[Section 13, proof of Theorem 2.8]{E24} to the one-phase setting.

Regarding the strong integrability assumption in Theorem \ref{main_theorem}, as discussed briefly in \cite{ERZ25}, in the context of minimal surface theory, White and Chan \cite{W89, C97} constructed complete minimal hypersurfaces asymptotic to minimizing non-quadratic cones. More precisely, given a strictly minimizing cone with an isolated singularity, Chan constructed a family of complete area-minimizing hypersurfaces whose parameters correspond to the slow-decaying Jacobi fields of the cone. These results suggest that analogous families may exist in the one-phase free boundary problem.

\subsection{Outline of the proof}
Let $(u_k)$ be a sequence of $(c_1, c_2, c_3)$-non-degenerate weak solutions converging to $U$. To prove Theorem \ref{main_theorem}, it suffices to show that its conclusion holds for $u_k$ for all sufficiently large $k$. By blowing up $(u_k)$ at a suitably chosen point and scale capturing the behavior near the singularity of $U$, we get a global non-degenerate weak solution $v_\infty$. Blowing-down $v_\infty$, we get a tangent cone $V$ at infinity of $v_\infty$ whose free boundary is a small graph over the free boundary of $U$, with graph function $\tilde{\xi}$. Once we know that $V$ is sufficiently close to $U$ in $L^2$, \cite[Corollary 4.7]{ERZ25} implies that
\[
V = U \circ  q
\]
for some $q \in O(d)$. The rigidity theorem \cite[Theorem 1.3]{ERZ25} then implies that, up to a rotation and a translation, $v_\infty$ is equal to one of $U, \overline{U}_t, \underline{U}_t$ for some $t > 0$. From this, we will be able to conclude that $u_k$ satisfies the conclusion of Theorem \ref{main_theorem} for all sufficiently large $k$.

The overall proof follows the strategy of \cite{E24}. The rigidity result \cite[Theorem 1.3]{ERZ25} provides the main tool needed to carry out this argument in the one-phase setting. However, there is one missing piece that we need to carefully justify. To apply \cite[Corollary 4.7]{ERZ25}, we need $\|U-V\|_{L^2}$ to be sufficiently small, whereas a priori, we only have small $\|\tilde{\xi}\|_{C^2}$. Therefore, we need a qualitative stability result showing that, if the graph of the free boundary of $V$ over the free boundary of $U$ is small, then their $L^2$ difference is also small. While the estimate 
\[
\|\tilde{\xi}\|_{C^2} \leq C\|U- V\|_{L^2}
\]
is readily available in \cite[Lemma 3.8, 3.10]{ERZ25}, the converse estimate requires a careful justification. This is proven in Section \ref{Graphical_stability_section}.

\subsection{AI Disclosure}
AI tools were used solely for language editing and manuscript review. The mathematical content and proof strategy were developed by the author.

\subsection{Acknowledgements}
The author thanks Luca Spolaor for suggesting this problem and Nicholas Edelen for a helpful discussion on simplifying the blow-up argument.

\section{Preliminaries}
\subsection{Non-degenerate weak solutions}
We start by recalling the definition of weak solution and non-degeneracy.

\begin{definition} \label{def_non_degenerate_weak_sol}
Let $\Omega \subseteq \mathbb R^d$ be open and $c_1, c_2, c_3 > 0$ with $c_1 \leq c_2$. We say $u : \Omega \rightarrow \mathbb R$ is a $(c_1, c_2, c_3)$-non-degenerate weak solution (of the one-phase Bernoulli problem) in $\Omega$ if
\begin{enumerate}
    \item \label{def_weak_sol_1} $u \in C(\Omega) \cap H^1_{\mathrm{loc}}(\Omega)$, $u \geq 0$ in $\Omega$ and $\Delta u = 0$ in $\{u > 0\}$.
    
    \item \label{def_weak_sol_2} $\{u > 0\}$ is a set of locally finite perimeter in $\Omega$ and $\Delta u = \mathcal{H}^{d-1}\mid_{\partial^*\{u > 0\}}$ in the sense of distributions, i.e.
    \[
    - \int_\Omega \nabla u \cdot \nabla \phi = \int_{\partial^*\{u > 0\}} \phi \, d\mathcal{H}^{d-1}
    \]
    for any $\phi \in C_c^\infty(\Omega)$.
    
    \item \label{def_non_deg} for any $x_0 \in \partial \{u > 0\}$, $r > 0$ with $B_r(x_0) \Subset \Omega$
    \[
        c_1 \leq \frac{1}{r} \intavg_{\partial B_r(x_0)} u \, d \mathcal{H}^{d-1} \leq c_2.
    \]
    
    \item \label{def_abundant_zero} for any $x_0 \in \partial \{u > 0\}$, $r > 0$ with $B_r(x_0) \Subset \Omega$
    \[
        |\{u = 0\} \cap B_r(x_0)| \geq c_3 |B_r|.
    \]
\end{enumerate}
We say $u$ is a $(c_1, c_2)$-non-degenerate weak solution if (\ref{def_weak_sol_1}), (\ref{def_weak_sol_2}), and (\ref{def_non_deg}) holds with $c_1, c_2$, and $u$ is a weak solution if (\ref{def_weak_sol_1}) and (\ref{def_weak_sol_2}) holds.
\end{definition}

Thorughout the paper, we will assume that $c_1, c_2, c_3 > 0$ with $c_1 \leq c_2$.

\begin{remark}
One can easily verify that if $u$ is a $(c_1, c_2, c_3)$-non-degenerate weak solution in $\Omega$, then for any open $U \subseteq \Omega$, $x_0 \in \mathbb R^d$ and $r > 0$, $u_{x_0, r}$ is a $(c_1, c_2, c_3)$-non-degenerate weak solution in $\frac{1}{r} (U - x_0)$.
\end{remark}

\begin{remark}
By \cite[Theorem 2.1 and Remark 2.3]{LW17}, if $u$ satisfies (\ref{def_weak_sol_1}) and (\ref{def_non_deg}) from Definition \ref{def_non_degenerate_weak_sol}, then $\{u > 0\}$ is a set of locally finite perimeter in $\Omega$. Hence, the first part of (\ref{def_weak_sol_2}) from Definition \ref{def_non_degenerate_weak_sol} is redundant.
\end{remark}

\begin{remark}
By \cite[Section 3]{AC81}, minimizers are $(c_1, c_2, c_3)$-non-degenerate weak solutions, where $c_1, c_2, c_3$ are dimensional constants. 
\end{remark}

\subsection{Uniform Lipschitz bounds and compactness}

We state and give proofs for uniform Lipchitz bounds and compactness results that will be fundamental for proving our result. The proofs are originally from \cite{AC81}.

\begin{lemma} \label{Lemma_uniform_lipschitz}
Let $\Omega$ be an open subset of $\mathbb R^d$. Suppose $u : \Omega \rightarrow [0, \infty)$ is a continuous function which is harmonic in $\Omega \cap \{u > 0\}$. Furthermore, suppose there exists $c_2 > 0$ such that for $x \in \partial\{u > 0\} \cap \Omega$ and $r > 0$ so that $B_r(x) \Subset \Omega$,
\[
\frac{1}{r} \intavg_{\partial B_{r}(x)} u \, d\mathcal{H}^{d-1} \leq c_2
\]
Then $u \in C^{0, 1}_{\mathrm{loc}}(\Omega)$.\\
Moreover, for a connected open subset $D \Subset \Omega$ such that $D \cap \partial \{u > 0\} \ne \emptyset$,
\[
\|u\|_{L^\infty(D)}, \mathrm{Lip}(u\mid_D) \leq C c_2
\]
where $C = C(d, D, \Omega) > 0$ which is independent of $u$.
\end{lemma}

\begin{proof}
Let $D \Subset D' \Subset \Omega$ where $D, D'$ are open. Let $r_0 := \frac{1}{4}\mathrm{dist}(D, \partial D') > 0$ and $C_0 := \frac{\sup_{D'} u}{r_0}$. For $x \in D \cap \{u > 0\}$, let $B_r(x)$ be the largest ball centered at $x$ which is contained in $D'\cap \{u > 0\}$. When $r \geq r_0$,
\[
\sup_{\partial B_r(x)} u \leq C_0 r
\]
Consider $r \leq r_0$. Then there exists $x_0 \in \partial B_r(x) \cap \partial \{u > 0\} \cap D'$. Since $B_{3r}(x_0) \subseteq D'$ and $x_0 \in \partial\{u> 0\} \cap \Omega$,
\[
\intavg_{\partial B_{3r}(x_0)} u \, d\mathcal{H}^{d-1} \leq 3c_2 r
\]
Fix $y \in \partial B_r(x)$. Since $u$ is subharmonic on $B_{3r}(x_0)$, $u \geq 0$ and $|y-z| \geq r$ for any $z \in \partial B_{3r}(x_0)$,
\[
u(y) \leq \int_{\partial B_{3r}(x_0)} \frac{(3r)^2 - |y-x_0|^2}{3r\omega_{d-1} |y - z|^d} u(z) \, d\mathcal{H}^{d-1}(z) \leq C_d \intavg_{\partial B_{3r}(x_0)} u \, d\mathcal{H}^{d-1}
\]
where $\omega_{d-1} = \mathcal{H}^{d-1}(\partial B_1)= d \, \mathcal{H}^d(B_1)$.\\
Therefore,
\[
\sup_{\partial B_{r}(x)} u \leq C_d c_2 r
\]
where $C_d > 0$ is a dimensional constant.\\
Hence, by the gradient estimate for harmonic function and the maximum principle,
\[
|\nabla u(x)| \leq \frac{C}{r^{d+1}} \int_{B_r(x)} u \leq C\max\{C_0, c_2\}
\]
where $C = C(d) > 0$.\\
Since $\nabla u = 0$ a.e. on $\{u = 0\}$, we conclude that $u \in W^{1, \infty}(D)$ with $\|\nabla u\|_{L^\infty(D)} \leq C\max\{C_0, c_2\}$.\\

\noindent
Now, fix a connected open subset $D \Subset \Omega$ so that $D \cap \partial \{u > 0\} \ne \emptyset$. Take open subsets $D \Subset D' \Subset D'' \Subset \Omega$, where we may take $D'$ to be connected as $D$ is connected. By the previous argument, it suffices to get the desired bound on $\|u\|_{L^\infty(D')}$.\\
Let $r_0 := \frac{1}{16}\mathrm{dist}(D', \partial D'') > 0$. Note that $r_0 = r_0(D, \Omega)$. Pick $x_0 \in D \cap \partial \{u > 0\}$. Fix an arbitrary $x_* \in D'$. Since $x_0, x_* \in D'$ and $D'$ is connected, there is a finite sequence $x_0, ..., x_L \in D'$ such that $x_L = x_*$, and $x_{j+1} \in B_{r_0}(x_{j})$ for $j = 0,...,L-1$ where $L \leq N(D', r_0) = N(D, \Omega)$. Let $j_0$ be the maximum over $\{j \in [0..L] : B_{2r_0}(x_j) \cap \partial \{u > 0\} \ne \emptyset \}$. Since $x_0 \in \partial \{u > 0\}$, $j_0$ is well-defined. Then for $j \geq j_0 + 1$, $u$ is harmonic on $B_{2r_0}(x_j)$ and $x_{j-1}, x_j \in B_{r_0}(x_j)$, by the Harnack's inequality,
\[
u(x_j) \leq C(d, r_0)u(x_{j-1})
\]
For $j = j_0$, pick $y_0 \in B_{2r_0}(x_{j_0}) \cap \partial \{u>0\}$. As above, since $u$ is subharmonic in $B_{8r_0}(y_0) \Subset D''$, $u \geq 0$, and $|x_{j_0} - z| \geq r_0$ for any $z \in \partial B_{8r_0}(y_0)$,
\[
u(x_{j_0}) \leq C_d \intavg_{\partial B_{8r_0}(y_0)} u\, d\mathcal{H}^{d-1} \leq C(d, r_0)c_2
\]
Combining the previous two inequalities,
\[
u(x_*) = u(x_L) \leq C(d, D, \Omega) c_2
\]
\end{proof}

The following lemma, in particular, shows that free boundary is a Lebesgue null set.

\begin{lemma}
Let $\Omega \subseteq \mathbb R^d$ be open. Suppose $u \in C(\Omega)$, $u \geq 0$ on $\Omega$, $u$ is harmonic on $\{u > 0\} \cap \Omega$, and $u$ satisfies the $(c_1, c_2)$ non-degeneracy condition from (\ref{def_non_deg}). Then $\mathcal{H}^{d-1}(\partial \{ u > 0\} \cap D) < \infty$ for any $D \Subset \Omega$.
\end{lemma}
\begin{proof}
\cite[Theorem 4.3, Theorem 4.5]{AC81}
\end{proof}

The following is the compactness lemma.

\begin{lemma}\label{Lemma_compactness} Let $\Omega$ be an open subset of $\mathbb R^d$. Suppose $u_i : \Omega \rightarrow [0, \infty)$ is a sequence of $(c_1, c_2)$-non-degenerate weak solutions. Also, suppose that $D \Subset \Omega$ is a connected open subset  satisfying
\[
\partial \{u_i > 0\} \cap D \ne \emptyset
\]
for all $i$.\\
Then there exists $u_0 : D \rightarrow [0, \infty)$ and passing to a subsequence
\begin{enumerate}
    \item $u_i \rightarrow u_0$ in $C^{0, \alpha}(D)$ for every $\alpha \in [0, 1)$

    \item $u_0$ is harmonic in $\{u_0 > 0\} \cap D$ with $u_0 \in W^{1, \infty}(D)$
    
    \item $\nabla u_i \rightarrow \nabla u_0$ weak-* in $L^\infty(D)$

    \item $\partial \{u_i > 0\} \rightarrow \partial \{u_0 > 0\}$ locally in Hausdorff distance in $D$

    \item $u_i \rightarrow u_0$ in $H^1(D)$
\end{enumerate}
\end{lemma}
\begin{proof}
Let $D', D''$ be a connected open set with $D \Subset D' \Subset D'' \Subset \Omega$.
By Lemma \ref{Lemma_uniform_lipschitz}, (1), (2), (3) follows immediately with $D$ replaced by $D''$.\\
We prove (4) with $D$ replaced by $D''$. For this, it suffices to show
\begin{enumerate}[label=(\alph*)]
    \item for any $B_r(x) \Subset D''$ disjoint from $\partial\{u_0 > 0\}$, $B_{r/2}(x)$ is eventually disjoint from $\partial\{u_i > 0\}$.
    \item for any $B_r(x) \Subset D''$ and any subsequence $i_j$, if $B_r(x)$ is disjoint from $\partial\{u_{i_j} > 0\}$ for all $j$, then $x \notin \partial\{u > 0\}$.
\end{enumerate}
Consider (a). If $B_r(x) \Subset D''$ is disjoint from $\partial\{u_0 > 0\}$, then either $u_0 > C$ on $B_{3r/4}(x)$ for some $C > 0$ or $u_0 \equiv 0$ on $B_{3r/4}(x)$. By the uniform convergence, $\partial \{u_i > 0\}$ is eventually disjoint from $B_{r/2}(x)$ for the first case. For the second case, suppose there are infinitely many $i$ so that $B_{r/2}(x)$ intersects $\partial \{u_i > 0\}$. Passing to a subsequence, pick $x_i \in B_{r/2}(x) \cap \partial \{u_i > 0\}$. By compactness, we may assume $x_i \rightarrow x_0$ for some $x_0 \in \overline{B}_{r/2}(x)$. Then
\[
\intavg_{\partial B_{s}(x_i)} u_i \, d\mathcal{H}^{d-1} \geq c_1 s
\]
for any $s \in (0, r/4)$.\\
By integrating with respect to $s$ from $0$ to $r/4$,
\[
\int_{B_{r/4}(x
_i)} u_i \geq C(d, r) c_1
\]
Taking limit as $i \rightarrow \infty$, we get a contradiction.\\
Now, consider (b). Since $B_r(x)$ is disjoint from $\partial \{u_{i_j} > 0\}$, $u_{i_j}$ is harmonic in $B_r(x)$. Then $u_0$ is harmonic in $B_r(x)$ which implies that $x \notin \partial \{u > 0\}$.\\

\noindent
Lastly, we prove the local $H^1$ convergence. By the local Hausdorff convergence of the free boundary,
\[
\partial \{u_i > 0\} \cap \overline{D'} \subseteq N_\epsilon := T_\epsilon(\partial \{u_0 > 0\} \cap D'') 
\]
for $i \gg 1$.\\
Hence, $u_i - u_0$ is harmonic in $D' \setminus \overline{N_\epsilon}$ for $i \gg 1$ and $D \setminus \overline{N_{2\epsilon}} \Subset D' \setminus \overline{N_\epsilon}$. By the Caccioppoli inequality
\[
\int_{D \setminus \overline{N_{2\epsilon}}} |\nabla (u_i - u_0) |^2 \leq C \int_{D' \setminus \overline{N_\epsilon}} |u_i - u_0|^2
\]
where $C = C(\delta) > 0$ where $\delta := d(D \setminus \overline{N_{2\epsilon}}, \partial(D' \setminus \overline{N_\epsilon}))$.\\
Also,
\[
\int_{D \cap \overline{N_{2\epsilon}}} |\nabla (u_i - u_0) |^2 \leq M \, \mathcal{H}^d(D\cap \overline{N_{2\epsilon}})
\]
where $M := \sup\{\|\nabla u_i\|_{L^\infty(D)} : i \geq 0\}$ which is finite by Lemma \ref{Lemma_uniform_lipschitz}.\\
Since $\mathcal{H}^d(\partial\{u_0 > 0\}) = 0$, we have $\mathcal{H}^d(D\cap \overline{N_{2\epsilon}}) \rightarrow 0$ as $\epsilon \rightarrow 0$. Therefore, $u_i \rightarrow u_0$ in $H^1(D)$.
\end{proof}

\subsection{$L^1$ convergence of the positivity sets}
We also need $L^1$ convergence of positivity set, i.e. $\chi_{\{u_l > 0\}} \rightarrow \chi_{\{u_0 > 0\}}$ in $L^1$. By \cite[Lemma 3.4]{AC81}, if $(u_i)$ a sequence of minimizers converging to $u_0$ locally uniformly, then $\chi_{\{u_i > 0\}} \rightarrow \chi_{\{u_0 > 0\}}$ in $L^1$. However, this is not true for a sequence of $(c_1, c_2, c_3)$-non-degenerate weak solutions as the following example shows
\begin{example}
Let $x_i := \left( 1 - \frac{2}{i} \right) e_d$, where $e_d = (0, 0, ..., 1)$.\\
Define
\[
u_i(x) =
\begin{cases} 
      0 & x \in B_{\frac{1}{i}}(x_i)\\
      \frac{1}{i(d-2)} \left( 1 - \left( \frac{1}{i|x-x_i|}\right)^{d-2} \right) & x \in B_1\setminus B_{\frac{1}{i}}(x_i)
\end{cases}
\]
One can easily verify that $(u_i)$ is a sequence of $(c_1, c_2, c_3)$-non-degenerate weak solutions converging uniformly to 0. Indeed, by symmetry, to verify (\ref{def_non_deg}), it suffices to show that
\[
\frac{1}{r}\sup_{B_r(x_i-\frac{1}{i}  e_d)} u_i =  \frac{1}{ir(d-2)} \left( 1 - \left( \frac{1}{1 + ir} \right)^{d-2} \right)
\]
is uniformly bounded below for $r \in (0, \frac{3}{k})$ and $i \in \mathbb N$.\\
Likewise, to show (\ref{def_abundant_zero}) it suffices to show that
\[
\frac{|B_{1/i} \cap B_r(\frac{1}{i}e_d)|}{|B_r|}
\]
is uniformly bounded below for $r \in (0, \frac{3}{k})$ and $i \in \mathbb N$.\\
By rescaling, this is equivalent to showing that
\[
\frac{|B_1 \cap B_r(e_d)|}{|B_r|}
\]
is uniformly bounded below for $r \in (0, 3)$.\\
\end{example}

Fortunately, the only issue is that the limit vanishing everywhere. The proof follows immediately from Kriventsov and Weiss' recent result \cite{KW25}. We provide a few details with a proof of the $L^1$ convergence when the limiting function $u_0$ is not identically equal to 0.

\begin{definition}
Let $\Omega \subseteq \mathbb R^d$ be open. We say a pair $(u, \chi)$ with $u : \Omega \rightarrow [0, \infty)$ and $\chi : \Omega \rightarrow \{0, 1\}$ is a variational solution in $\Omega$ with constant $C_V$ if
\begin{enumerate}
    \item $u \in C^{0, 1}(\Omega) \cap C^2(\{u > 0\})$ with $|\nabla u|\leq C_V$ in $\Omega$.
    \item $\chi$ is Borel measurable
    \item $\chi_{\{u > 0\}}\leq \chi$ $\mathcal{L}^{d}$-a.e. in $\Omega$.
    \item For every vector field $\phi \in C_c^{0, 1}(\Omega ; \mathbb R^d)$,
    \[
        \int_{\Omega} (|\nabla u|^2 + \chi)\mathrm{div}\phi - 2\nabla u \cdot D\phi \nabla u = 0.
    \]
\end{enumerate}
\end{definition}

\begin{lemma} \label{Lemma_positivity_set_convergence}
Let $\Omega$ be an open subset of $\mathbb R^d$. Suppose $u_i : \Omega \rightarrow [0, \infty)$ is a sequence of $(c_1, c_2)$-non-degenerate weak solutions and $u_0 : \Omega \rightarrow [0, \infty)$ be a (continuous) function so that $u_i \rightarrow u_0$ in $C_{\mathrm{loc}}(\Omega) \cap H_{\mathrm{loc}}^1(\Omega)$. Let $D \Subset \Omega$ be a connected open subset such that $u_0 \not\equiv 0$ in $D$. Then
\[
    \chi_{\{u_i > 0\}} \rightarrow \chi_{\{u_0 > 0\}}
\]
in $L^1(D)$.\\
Moreover, $u_0$ satisfies the $(c_1, c_2)$ non-degeneracy condition, and $(u_0, \chi_{\{u_0 > 0\}})$ is a variational solution in $D$. If, in addition, $u_i$ is a sequence of $(c_1, c_2, c_3)$ non-degenerate weak solutions, then $u_0$ is a $(c_1, c_2, c_3)$ non-degenerate weak solution in $D$.
\end{lemma}
\begin{proof}
Take an open connected set $D'$ such that $D \Subset D' \Subset \Omega$. By the proof of \cite[Theorem 5.1]{W98}, $(u_i, \chi_{\{u_i > 0\}})$ is a variational solution in $D'$ for all $i$ with the same constant $C_V$ due to Lemma \ref{Lemma_uniform_lipschitz}.

By \cite[Lemma 3.3]{KW25}, $\chi_{\{u_i > 0\}}$ is in $BV(D)$ with $\|\chi_{\{u_i > 0\}}\|_{BV(D)} \leq C$ where $C > 0$ is independent of $i$. Hence, by the compactness of BV functions, $\chi_{\{u_i > 0\}} \rightarrow \chi$ in $L^1(D)$ for some $\chi \in BV(D)$. In particular, by considering a subsequence converging pointwise, $\chi \in \{0, 1\}$ and $\chi_{\{u > 0\}} \leq \chi$ $\mathcal{L}^d$-a.e. in $D$. Also, by the $H^1(D)$ convergence of $(u_i)$, $L^1(D)$ convergence of $(\chi_i)$ and uniform bound on $|\nabla u|$ in $D$, we conclude that $(u, \chi)$ is a variational solution in $D$.

By \cite[Theorem 8.6]{KW25}, we conclude that $\chi = \chi_{\{u > 0\}}$ $\mathcal{L}^{d}$-a.e. The $(c_1, c_2)$-non-degeneracy follows from taking limit. Also, if we additionally assume that $u_i$ is a sequence of $(c_1, c_2, c_3)$ non-degenerate weak solutions, by the proof of \cite[Theorem 5.1]{W98} (also, \cite[Remark 3.2]{ERZ25}), we conclude that $u_0$ is a weak solution. 
\end{proof}

\subsection{Weiss' boundary-adjusted energy}
Recall that for $u \in H_{\mathrm{loc}}^1(\Omega)$ and $x \in \Omega$, Weiss' boundary adjusted energy at $x$ with scale $r \in (0, \mathrm{dist}(x, \partial \Omega))$ is defined as
\[
    W_u(x, r) := \frac{1}{r^d}\int_{B_r(x)} |\nabla u|^2 \, dx - \frac{1}{r^{d+1}} \int_{\partial B_r(x)} u^2 \, d\mathcal{H}^{d-1} + \frac{1}{r^d} |\{u > 0\} \cap B_r(x)|.
\]
For convenience, we write $W_u(x) := \lim_{r \downarrow 0} W_u(x, r)$ which is well-defined by the following monotonicity formula.
\begin{lemma} \label{Lemma_monotonicity}
Let $u : \Omega \rightarrow [0, \infty)$ be a $(c_1, c_2)$-non-degenerate weak solution. Then for $x_0 \in \Omega$ and $r \in (0, \mathrm{dist}(x_0, \partial\Omega))$,
\[
    \frac{d}{dr} W_u(x_0,r) = \frac{2}{r^{d+2}} \int_{\partial B_r(x_0)} ((x-x_0) \cdot \nabla u(x) - u(x))^2\, d\mathcal{H}^{d-1}(x).
\]
In particular, $W_u(x, \cdot)$ is monotone increasing and $W_u(x,\cdot)$ is constant on $(r_1, r_2)$ if and only if $u(\cdot + x)$ is one-homogeneous on $A_{r_1, r_2}$.
\end{lemma}
\begin{proof}
\cite[Theorem 3.1, Theorem 5.1]{W98}
\end{proof}

\begin{corollary} \label{Corollary_Weiss_controls_homogeneity}
Suppose $u: B_1 \rightarrow [0, \infty)$ is a $(c_1, c_2)$-non-degenerate weak solution. 
For $r_1, r_2, \rho > 0$ with $r_1 < r_2$ so that $\rho r_2 < 1$,
\[
    \int_{A_{r_1, r_2}} ( x\cdot \nabla u_\rho(x) - u_\rho(x) ) ^2\, dx \leq \frac{r_2^{d+2}}{2} ( W_u(0, \rho r_2) - W_u(0, \rho r_1) )
\]
Moreover, for $0 < s < r < 1$,
\[
    \int_{\partial B_1} |u_r(x) - u_s(x)|^2\, d\mathcal{H}^{d-1}(x) \leq \frac{r-s}{2s} (W_u(0, r) - W_u(0, s))
\]
\end{corollary}
\begin{proof}
Wlog, assume $\rho = 1$.
\begin{equation}
\begin{split}
    \int_{A_{r_1, r_2}} ( x\cdot \nabla u(x) - u(x) ) ^2\, dx
    &= \int_{r_1}^{r_2} \int_{\partial B_s} (x \cdot \nabla u (x) - u(x))^2 \, d\mathcal{H}^{d-1}(x) \, ds\\
    &= \int_{r_1}^{r_2} \frac{s^{d+2}}{2} \frac{d}{ds} W_u(0, s)\, ds \\
    &\leq \frac{r_2^{d+2}}{2} ( W_u(0, r_2) - W_u(0, r_1) )
\end{split}
\end{equation}

For the second estimate, for $x \in \partial B_1$,
\[
    \frac{d}{dt} u_t(x) = \frac{tx\cdot \nabla u (tx) - u(tx)}{t^2}
\]
for a.e. $t \in (s, r)$.

Then
\[
    |u_r(x) - u_s(x)|^2 \leq (r-s) \int_s^r \frac{|tx\cdot \nabla u(tx) - u(tx)|^2}{t^4}\, dt
\]
and so
\[
    \int_{\partial B_1} |u_r(x) - u_s(x)|^2\, d\mathcal{H}^{d-1}(x) \leq (r-s) \int_s^r \frac{1}{t^4} \int_{\partial B_1} |tx\cdot \nabla u(tx) - u(tx)|^2 \, dx\, dt
\]
By Lemma \ref{Lemma_monotonicity},
\[
    \int_{\partial B_1} |u_r(x) - u_s(x)|^2\, d\mathcal{H}^{d-1}(x) \leq \frac{r-s}{2}\int_s^r \frac{1}{t} \frac{d}{dt} W_u(0, t) \, dt \leq \frac{r-s}{2s} (W_u(0, r) - W_u(0, s))
\]
\end{proof}

The following characterizes the $W_U(x)$ when $U$ is a minimizing cone with $\mathrm{Sing}(U) = \{0\}$. To begin with, we recall the definition of a minimizing cone.

\begin{definition}
We say $U : \mathbb R^d \rightarrow [0, \infty)$ is a minimizing cone if $U \in H^1_{\mathrm{loc}}(\mathbb R^d)$ is a minimizer of Alt-Caffarelli energy $\mathcal{F}$ i.e. $\mathcal{F}(U; B_R) \leq \mathcal{F}(v; B_R)$ for any $R > 0$ and $v \in H^1(B_R)$ with $v - U \in H^1_0(B_R)$ and $U$ is one-homogeneous (i.e. $U(rx) = rU(x)$ for any $x \in \mathbb R^d$ and $r > 0$).
\end{definition}

\begin{lemma} \label{Lemma_characterization_of_Weiss_for_minimizer}
Let $U$ be a minimizing cone with $\mathrm{Sing}(U) = \{0\}$. Then for $x \in \mathbb R^d$, 
\[
    \text{$W_U(x) > \frac{\omega_d}{2}$ if and only if $x = 0$}
\]
where $\omega_d := \mathcal{H}^d(B_1)$.
\end{lemma}
\begin{proof}
By \cite[Lemma 9.20, Lemma 9.22]{V23}, $W_U(0) > \frac{\omega_d}{2}$ and for $x \in \partial \{U > 0\}$, $W_U(x) > \frac{\omega_d}{2}$ implies that $x = 0$. From the definition of the Weiss energy, one can easily verify that $W_U(x) = 0$ if $x \in \{U = 0\}^\circ$, and $W_U(x) = -\infty$ if $x \in \{ U > 0\}$.
\end{proof}

\section{Graphical stability of one-homogeneous solutions} \label{Graphical_stability_section}
We define a tubular neighborhood of $\partial \{U > 0\}$.
\begin{definition} \label{Definition_tubular_neighborhood}
Let $U$ be a minimizing cone with $\mathrm{Sing}(U) = \{0\}$,  $\tau> 0$ and $0 \leq r_1 < r_2$, we define
\[
    C := \{\omega + t\nabla U(\omega) : \omega \in \partial \{U > 0\} \cap \mathbb S^{d-1}, |t| < \tau\}
\]
\[
    C_{r_1, r_2} := \{r\omega : \omega \in C, r_1 < r < r_2 \}
\]
\end{definition}

In particular, there exists $\tau > 0$ sufficiently small depending on $U$ such that $C_{0, \infty}$ is a tubular neighborhood of $\partial \{U > 0\} \setminus \{0\}$.

In \cite[Lemma 3.10]{ERZ25}, the authors prove a quantitative estimate showing that small $\|u - U\|_{C^2}$ implies small graphical distance between their free boundaries. In this section, we prove a partial converse: Theorem \ref{Theorem_SMALL_FREE_BOUNDARY_GRAPH_IMPLIES_SMALL_L^2} 
shows that if $U$ and $V$ are one-homogeneous weak solutions and the free boundary of $V$ is a sufficiently small $C^2$ graph over the free boundary of $U$, then $\|U - V\|_{L^2}$ is small.

\begin{lemma} \label{SPHERICAL_GRAPH}
There exists $\epsilon_0(U), C(U)>0$ so that the following holds: Suppose $\xi \in C^{2, \alpha}(\partial \{U > 0\} \cap A_{\frac{1}{2}, \frac{3}{2}})$ and $u : A_{\frac{1}{2}, \frac{3}{2}} \rightarrow [0, \infty)$ is continuous  so that
\[
\partial \{ u> 0\} \cap A_{\frac{1}{2}, \frac{3}{2}} = \mathrm{graph}_{\partial \{U > 0\} \cap A_{\frac{1}{2}, \frac{3}{2}}}(\xi) \cap A_{\frac{1}{2}, \frac{3}{2}}
\]
and for each $x_0 \in \partial \{ U > 0\} \cap A_{\frac{1}{2}, \frac{3}{2}}$, there exists open neighborhoods $E, F$ of $x_0, y_0 := x_0 + \xi(x_0)\nu_{x_0}$, respectively, so that
\[
\{u > 0\} \cap F \cap A_{\frac{1}{2}, \frac{3}{2}} = \{x + t \nu_x : x\in \partial \{U>0\} \cap E \cap A_{\frac{1}{2}, \frac{3}{2}}, t > \xi(x)\} \cap F \cap A_{\frac{1}{2}, \frac{3}{2}}
\]
If $\|\xi\|_{C^1(\partial \{U > 0\} \cap A_{\frac{1}{2}, \frac{3}{2}})} \leq \epsilon_0$, then there exists $\eta \in C^{2,\alpha}(\partial \{U > 0\} \cap \mathbb S^{d-1})$ so that $\|\eta\|_{C^{2,\alpha}(\partial \{U > 0\} \cap \mathbb S^{d-1})} \leq C \|\xi\|_{C^{2,\alpha}(\partial \{U > 0\} \cap A_{\frac{1}{2}, \frac{3}{2}})}$ and
\[
\partial \{ u > 0\} \cap \mathbb S^{d-1} = \{\mathrm{exp}_{\theta}^{\mathbb S^{d-1}}(\eta(\theta)\nu_\theta) : \theta \in \partial \{U > 0\} \cap \mathbb S^{d-1}\}
\]
and for each $\theta_0 \in \partial \{U>0\} \cap \mathbb S^{d-1}$, there exists open neighborhoods $E, F$ of $\theta_0, \omega_0 := \exp_{\theta_0}^{\mathbb S^{d-1}}(\eta(\theta_0)\nu_{\theta_0})$, respectively, so that
\[
\{ u > 0 \} \cap F \cap \mathbb S^{d-1} = \{\mathrm{exp}_{\theta}^{\mathbb S^{d-1}}(\tau\nu_\theta) : \theta \in \partial \{U > 0\} \cap E \cap \mathbb S^{d-1}, \tau > \eta(\theta). |\tau| < \tau_*\} \cap F
\]
\end{lemma}
\begin{proof}
Let $F(r, \theta) := |r\theta + \xi(r\theta) \nu_{\theta}|^2 = r^2 + \xi(r\theta)^2$ for $r \in \left( \frac{1}{2},\frac{3}{2}\right), \theta \in \partial \{U>0\} \cap \mathbb S^{d-1}$, where we note that $\nu_{\theta} = \nu_{r\theta}$. Then
\[
\partial \{u > 0\} \cap \mathbb S^{d-1} = \{r\theta + \xi(r\theta) \nu_{\theta} : F(r, \theta) = 1\}
\]
Also, $\partial_r F(r, \theta) = 2r + 2\xi(r\theta) \partial_r \xi(r\theta)$. For sufficiently small $\epsilon_0$, we have $\partial_r F(r, \theta) > \frac{1}{2}$ for $r \in \left( \frac{1}{2},\frac{3}{2}\right), \theta \in \partial \{U>0\} \cap \mathbb S^{d-1}$ and $F(\frac{1}{2} +) < 1 < F(\frac{3}{2}-)$. Hence, for each $\theta \in \partial \{U>0\} \cap \mathbb S^{d-1}$, there is a unique $r(\theta) \in \left( \frac{1}{2},\frac{3}{2}\right)$ so that $F(r(\theta), \theta) = 1$. By the implicit function theorem, $r \in C^2(\partial \{U>0\} \cap \mathbb S^{d-1})$.\\
Now, $\mathrm{exp}_{\theta}^{\mathbb S^{d-1}}(\tau \nu_\theta) = \cos(\tau) \theta + \sin(\tau) \nu_\theta$. Hence, define $\eta(\theta) := \arcsin(\xi(r(\theta)\theta))$. Notice that $r(\theta) = \cos(\eta(\theta))$. By computations, the $C^2$ estimate and the $\eta$ spherical graphicality of the free boundary follows  by choosing sufficiently small $\epsilon_0, C > 0$ depending on $U$.\\

\noindent
Now, fix $\theta_0 \in \partial \{ U > 0\} \cap \mathbb S^{d-1}$ and let $\omega_0 := \exp_{\theta_0}^{\mathbb S^{d-1}}(\eta(\theta_0)\nu_{\theta_0})$. Then there exists $x_0 \in \partial \{ U > 0\} \cap A_{\frac{1}{2}, \frac{3}{2}}$ so that
\[
\omega_0 = x_0 + \xi(x_0) \nu_{x_0}
\]
By the assumption, we pick open neighborhoods $E, F$ of $x_0, \omega_0$, respectively, satisfying
\[
\{u > 0\} \cap F \cap A_{\frac{1}{2}, \frac{3}{2}} = \{x + t \nu_x : x\in \partial \{U>0\} \cap E \cap A_{\frac{1}{2}, \frac{3}{2}}, t > \xi(x)\} \cap F \cap A_{\frac{1}{2}, \frac{3}{2}}
\]
By shrinking $\epsilon_0$, we may assume that a tubular neighborhood $N$ of $\partial \{U > 0\} \cap \mathbb S^{d-1}$ contains $\partial \{ u > 0\} \cap \mathbb S^{d-1}$.\\
Let $F_1 := \{x + t \nu_x : x \in \partial \{ U > 0\} \cap E \cap A_{\frac{1}{2}, \frac{3}{2}}, t \in \mathbb R\}$. We define $F' := F_1 \cap F \cap N \cap \mathbb S^{d-1}$. Clearly, $F'$ is non-empty and open in $\mathbb S^{d-1}$. Let $E'$ be the projection of $F'$ onto $\partial\{U > 0\} \cap \mathbb S^{d-1}$. Then $E'$ is non-empty and open in $\partial \{U > 0\} \cap \mathbb S^{d-1}$. We show that
\[
\{ u > 0 \} \cap F' \cap \mathbb S^{d-1} = \{\mathrm{exp}_{\theta}^{\mathbb S^{d-1}}(\tau\nu_\theta) : \theta \in \partial \{U > 0\} \cap E' \cap \mathbb S^{d-1}, \tau > \eta(\theta), |\tau| < \tau_* \} \cap F'
\]
Fix $\omega_1 \in \{ u > 0 \} \cap F' \cap \mathbb S^{d-1}$. Then
\[
\omega_1 = x_1 + t_1 \nu_{x_1}
\]
for some $x_1 \in \partial \{ U> 0\} \cap E \cap A_{\frac{1}{2}, \frac{3}{2}}$ and $t_1 > \xi(x_1)$.\\
Since $\omega_1 \in N$, there is a unique $\theta_1 \in \partial \{ U > 0\} \cap \mathbb S^{d-1}$ and $|\tau_1|\leq C(\epsilon_0)$ such that
\[
\omega_1 = \exp_{\theta_1}(\tau_1 \nu_{\theta_1}) = \cos(\tau_1) \theta_1 + \sin(\tau_1) \nu_{\theta_1}
\]
where $C(\epsilon_0) \rightarrow 0$ as $\epsilon_0 \rightarrow 0$.\\
Hence, $x_1 = \cos(\tau_1)\theta_1$ and $t_1 = \sin(\tau_1)$. Let $g(\tau) := \sin(\tau) - \xi(\cos(\tau)\theta_1)$ defined for $|\tau| < C(\epsilon_0)$. We have $g(\tau_1) = t_1 - \xi(x_1) > 0$ and $g(\eta(\theta_1)) = \sin(\eta(\theta_1)) - \xi(\cos(\eta(\theta_1))\theta_1) = 0$. Now, $g'(\tau) = \cos(\tau) + O(\epsilon_0)$, so by shrinking $\epsilon_0$, if necessary, $g' > 0$ on $|\tau| < C(\epsilon_0)$, and so strictly increasing. Therefore, $\tau_1 > \eta(\theta_1)$.\\

\noindent
Now, let $\omega_1 \in \{\mathrm{exp}_{\theta}^{\mathbb S^{d-1}}(\tau\nu_\theta) : \theta \in \partial \{U > 0\} \cap E' \cap \mathbb S^{d-1}, \tau > \eta(\theta), |\tau| < \tau_*\} \cap F'$. Then
\[
\omega_1 = \cos(\tau_1) \theta_1 + \sin(\tau_1) \nu_{\theta_1} = x_1 + t_1 \nu_{x_1}
\]
for some $\theta_1 \in \partial \{U > 0\} \cap E' \cap \mathbb S^{d-1}$, $\tau_1 > \eta(\theta_1)$, and $x_1 \in \partial \{U > 0\} \cap E \cap A_{\frac{1}{2}, \frac{3}{2}}$, $t_1 \in \mathbb R$. It suffices to show that $t_1 > \xi(x_1)$. Using the same $g$, as before, we have $g(\eta(\theta_1)) = 0$ and since $g$ is strictly increasing and $\tau_1 > \eta(\theta_1)$, we have $g(\tau_1) > 0$, which is $t_1 > \xi(x_1)$.
\end{proof}

\begin{lemma} \label{DIFFEOMORPHISM_OF_POSITIVE_SET}
Let $\Gamma, \Lambda \subseteq \mathbb S^{d-1}$ be open with smooth boundary, and $\eta \in C^{2,\alpha}(\partial \Gamma)$, $\nu$ be the smooth inner unit normal on $\partial \Gamma$. There exists $\epsilon_0(\Gamma), C(\Gamma) > 0$ so that the following holds: Suppose $\|\eta\|_{C^0(\partial \Gamma)} \leq \epsilon_0$ and $\partial \Lambda$ is contained in a tubular neighborhood of $\partial\Gamma$ and  for each $x_0 \in \partial \Gamma$, there is an open neighborhood $E, F$ of $\theta_0, \omega_0 := \exp_{\theta_0}^{\mathbb S^{d-1}}(\eta(\theta_0)\nu_{\theta_0})$ so that 
\[
\partial \Lambda = \{\exp_\theta^{\mathbb S^{d-1}}(\eta(\theta)\nu_\theta) : \theta \in \partial \Gamma\}
\]
and
\[
\Lambda \cap F = \{\exp_\theta^{\mathbb S^{d-1}}(\tau\nu_\theta) : \theta \in \partial \Gamma \cap E, \tau > \eta(\theta), |\tau| < \tau_*\} \cap F
\]
Then there is a $C^{2,\alpha}$ diffeomorphism $\Phi : \overline{\Gamma} \rightarrow \overline{\Lambda}$ with $\|\exp_{\mathrm{Id}}^{-1}\Phi\|_{C^{2,\alpha}(\overline{\Gamma}; T\mathbb S^{d-1})} \leq C\|\eta\|_{C^{2,\alpha}(\partial \Gamma)}$
\end{lemma}
\begin{proof}
Since $\partial \Gamma$ is smooth and compact, let $N$ be a tubular neighborhood so that $\Psi(\theta, \tau) := \exp_{\theta}(\tau\nu_\theta)$ is a diffeomorphism from $\partial \Sigma \times (-\delta_0, \delta_0)$ to $N$.\\
Let $\chi : \mathbb R \rightarrow [0,1]$ be a smooth cut-off function so that $\chi(\tau) = 1$ for $|\tau| \leq \delta_0/2$, $\chi(\tau) = 0$ for $|\tau| \geq \delta_0$, and $|\chi'| \leq 1/\delta_0$. Also, set $\epsilon_0 := \frac{\delta_0}{4}$ so that $\|\eta\|_{C^0(\partial \Gamma)} \leq \frac{\delta_0}{4}$.\\
We define $\Phi : \overline\Gamma \rightarrow \mathbb S^{d-1}$ as follows. When $p \in N$, we have $(\theta, \tau) = \Psi^{-1}(p)$. 
Define 
\[
\Phi(p) := \exp_{\theta}^{\mathbb S^{d-1}}((\tau + \eta(\theta)\chi(\tau)) \nu_\theta)
\]
When $p \notin N$, we let $\Phi(p) := p$.\\
Since $\|\eta\|_{C^0(\partial \Gamma)} \leq \frac{\delta_0}{4}$, $g(\tau):=  \tau + \eta(\theta)\chi(\tau)$ is strictly increasing on $(-\delta_0, \delta_0)$. Also, $g((-\delta_0, \delta_0)) = (-\delta_0, \delta_0)$. Thus, $\Phi(p) \in N$ for $p \in N$. Hence, $\Phi$ is injective.\\

\noindent
To begin with, show that $\Gamma \setminus N = \Lambda \setminus N$. Suppose $p \in \Gamma \setminus N$. Let $B_{\tau_0}(p)$ be the geodesic ball in $\mathbb S^{d-1}$ such that $B_{\tau_0}(p) \subseteq \Gamma$ and $\partial B_{\tau_0}(p) \cap \partial \Gamma$ intersects at $\theta_0$. Then, $\tau_0 \geq \delta_0$ and $p = \exp_{\theta_0}^{\mathbb S^{d-1}}(\tau_0 \nu_{\theta_0})$ with $\exp_{\theta_0}^{\mathbb S^{d-1}}(\tau \nu_{\theta_0}) \in \Gamma$ for all $\tau \in (0, \tau_0]$. By the assumption, there is $\kappa_0 > 0$ such that $\exp_{\theta_0}^{\mathbb S^{d-1}}(\tau\nu_{\theta_0}) \in \Lambda$ for $\tau \in (\eta(\theta_0), \eta(\theta_0) + \kappa_0)$. Let $\tau_* := \max\{\tilde{\tau} > \eta(\theta_0) : \exp_{\theta_0}^{\mathbb S^{d-1}}(\tau\nu_{\theta_0}) \in \Lambda \text{ for $t\in (\eta(\theta_0), \tilde{\tau})$}\}$. Then $\tau_* \geq \eta(\theta_0) + \kappa_0$ and $\exp_{\theta_0}^{\mathbb S^{d-1}}(\tau_*\nu_{\theta_0}) \in \partial \Lambda$. Hence, $\tau_* \geq \delta_0$, otherwise, $\eta$ must be 2-valued, which is not the case. We would like to show that $\tau_* > \tau_0$. Suppose not. Then since $\|\eta\|_{C^0(\partial \Gamma)} \leq \frac{\delta_0}{4}$ and $\partial \Lambda$ is graphical over $\partial \Gamma$ by $\eta$, the geodesic ball $B_{\frac{\delta_0}{2}}(\exp_{\theta_0}^{\mathbb S^{d-1}}(\tau_*\nu_{\theta_0}))$ must intersect with $\partial \Gamma$. But since $\tau_* \in [\delta_0, \tau_0]$, we have $B_{\frac{\delta_0}{2}}(\exp_{\theta_0}^{\mathbb S^{d-1}}(\tau_*\nu_{\theta_0})) \subseteq B_{\tau_0}(p)$ which is contained in $\Gamma$, a contradiction. This proves $\Gamma \setminus N \subseteq \Lambda \setminus N$. With a similar argument, we get the equality.\\

\noindent
Hence, we are left with proving $\Phi(\Gamma \cap N) = \Lambda \cap N$. Let $\theta_0 \in \partial \Gamma$ and $\tau_0 \in (0, \delta_0)$. Since $g(\tau) = \tau + \eta(\theta_0)\chi(\tau)$ is strictly increasing and $g(0) = \eta(\theta_0)$, we have $\eta(\theta_0) < \tau_0 + \eta(\theta_0)\chi(\tau_0)$. By the assumption, there is $\kappa_0 > 0$ such that $\exp_{\theta_0}^{\mathbb S^{d-1}}(\tau \nu_{\theta_0}) \in \Lambda$ for $\tau \in (\eta(\theta_0), \eta(\theta_0)+ \kappa_0)$. By the same reasoning as above, $\exp_{\theta_0}^{\mathbb S^{d-1}}(\tau \nu_{\theta_0}) \in \Lambda$ for $\tau \in (\eta(\theta_0), \delta_0)$ (otherwise, $\eta$ is 2-valued), and so $\exp_{\theta_0}^{\mathbb S^{d-1}}((\tau_0 + \eta(\theta_0)\chi(\tau_0)) \nu_{\theta_0}) \in \Lambda$, proving that $\Phi(\Gamma \cap N) \subseteq \Lambda \cap N$.\\
Now, suppose $p \in \Lambda \cap N$. Then $p = \exp_{\theta_0}^{\mathbb S^{d-1}}(\tau_0 \nu_{\theta_0})$ for some $\theta_0 \in \partial \Gamma$ and $\tau_0 \in (-\delta_0, \delta_0)$. In particular, there exists $\tau_1 \in (-\delta_0, \delta_0)$ such that $\tau_0 = g(\tau_1) = \tau_1 + \eta(\theta_0)\chi(\tau_1)$. We need to show that $\tau_1 > 0$ so that $\Psi(\theta_0, \tau_1) \in \Gamma$. Suppose not. Since $g$ is increasing and $g(0) = \eta(\theta_0)$, we have $-\delta_0 < \tau_1 + \eta(\theta_0)\chi(\tau_1) \leq \eta(\theta_0)$. Now, by the assumption and the graphicality, $\exp_{\theta_0}^{\mathbb S^{d-1}}(\tau \nu_{\theta_0}) \notin \Lambda$ for $\tau \in (-\delta, \eta(\theta_0)]$. In particular, $p = \exp_{\theta_0}^{\mathbb S^{d-1}}((\tau_1 + \eta(\theta_0)\chi(\tau_1)) \nu_{\theta_0}) \notin \Lambda$, a contradiction.\\
Therefore, $\Phi : \overline{\Gamma} \rightarrow \overline{\Lambda}$ is a diffeomorphism. By computation, $\|\exp_{\mathrm{Id}}^{-1}\Phi\|_{C^2(\overline{\Gamma}; T\mathbb S^{d-1})} \leq C\|\eta\|_{C^2(\partial \Gamma)}$ follows.
\end{proof}

\begin{theorem} \label{Theorem_SMALL_FREE_BOUNDARY_GRAPH_IMPLIES_SMALL_L^2}
Suppose $U$ is a 1-homogeneous weak solutions where $\mathrm{Sing}(U) = \{0\}$. For any $\epsilon > 0$, there is $\delta(U, \epsilon) > 0$ so that the following holds: Suppose $V$ is a 1-homogeneous weak solution and $\xi \in C^{2,\alpha}(\partial\{U > 0\}\cap A_{1,2})$ with $\|\xi\|_{C^{2,\alpha}(\partial\{U > 0\}\cap A_{1,2})}< \delta$
\[
\partial\{ V > 0 \} \cap A_{1,2} = \mathrm{graph}_{\partial\{U > 0\}\cap A_{1,2}}(\xi) \cap A_{1,2}
\]
and for each $x_0 \in \partial \{ U > 0 \} \cap A_{1,2}$, there is an open neighborhood $E$ of $x_0$ so that
\[
\{ V > 0 \} \cap A_{1,2} \cap E = \{x + t\nu_x : x \in \partial \{U> 0\} \cap A_{1,2}, t > \xi(x) \} \cap E
\]
Then $\|U - V\|_{L^2(A_{1, 2})} < \epsilon$.
\end{theorem}
\begin{proof}
Suppose not. Then we have a sequence of 1-homogeneous weak solutions $V_i$ and $\xi_i \in C^{2,\alpha}(\partial\{U> 0\} \cap A_{1,2})$ with 
$\|\xi_i\|_{C^{2,\alpha}(\partial\{U> 0\} \cap A_{1,2})} \rightarrow 0$. Let $\Gamma := \{U > 0\} \cap \mathbb S^{d-1}$, $\Lambda_i := \{V_i > 0\} \cap \mathbb S^{d-1}$, $\partial \Gamma := \partial \{U > 0\} \cap \mathbb S^{d-1}$ and $\partial \Lambda_i := \partial \{V > 0\} \cap \mathbb S^{d-1}$. By Lemma \ref{SPHERICAL_GRAPH} and Lemma \ref{DIFFEOMORPHISM_OF_POSITIVE_SET}, there exists $\eta_i \in C^{2,\alpha}(\partial \Gamma)$ and a $C^{2,\alpha}$ diffeomorphism $\Phi_i : \overline\Gamma \rightarrow \overline \Lambda_i$ such that $\partial \Lambda$ is a spherical graph of $\eta_i$ over $\partial \Gamma$ and $\|\exp_{\mathrm{Id}}^{-1}\Phi_i\|_{C^{2,\alpha}(\overline\Gamma)} \leq \|\eta_i\|_{C^{2,\alpha}(\partial \Gamma)}$ with $\|\eta_i\|_{C^{2,\alpha}(\partial \Gamma)} \rightarrow 0$.\\
Let $v_i := V_i \mid_{\overline\Lambda_i}$. Since $V_i$ is 1-homogeneous weak solution,
\[
\Delta_{\mathbb S^{d-1}} v_i = -(d-1)v_i \quad\text{on $\Lambda_i$}
\]
\[
v_i = 0 \quad\text{on $\partial\Lambda_i$}
\]
\[
\partial_{\nu_i} v_i = 1 \quad\text{on $\partial\Lambda_i$}
\]
where $\nu_i$ is the unit inner normal of $\partial\Lambda_i$.\\
By differentiating $V_i(r\theta) = rV_i(\theta)$ with respect to $r$, we get $V_i(\theta) = \theta \cdot \nabla V_i(r\theta)$ for any $\theta \in \Lambda_i$. Since $V_i$ is harmonic in $\{V_i > 0\}$ and one-homogeneous, we have $\Delta_{\mathbb S^{d-1}} |\nabla V_i|^2 = 2|\nabla^2 V_i|^2$. Since $|\nabla V_i|^2$ is subharmonic in $\Lambda_i$ and $|\nabla V_i|^2 = 1$ on $\partial \Lambda_i$, we conclude that $|\nabla V_i|^2 \leq 1$ in $\overline{\Lambda_i}$. In particular,
\begin{equation}
    |v_i(\theta)| = |\theta \cdot \nabla V_i(\theta)| \leq 1
\end{equation}
for $\theta \in \overline{\Lambda_i}$. Together with the boundary Schauder estimate, we conclude that $\sup_i \|v_i\|_{C^{2,\alpha}(\overline{\Lambda_i})} < \infty$.

Define $w_i : \overline{\Gamma} \rightarrow [0, \infty)$ as $w_i := v_i \circ \Phi_i$. From the definition of $\Phi_i$, we have $\Delta_{\mathbb S^{d-1}} w_i  (x) = \Delta_{\mathbb S^{d-1}} v_i (\Phi_i(x))$ for $x$ away from the tubular neighborhood. For $x$ in the tubular neighborhood, by computing on the normal chart given by the tubular neighborhood,
\[
\Delta_{\mathbb S^{d-1}} w_i  = \Delta_{\mathbb S^{d-1}} v_i (\Phi_i) + R_i
\]
where
\[
    \|R_i\|_{C^0(\overline{\Gamma})} \leq C \|\eta_i\|_{C^2} \|v_i\|_{C^2(\overline{\Lambda_i})}.
\]
Likewise, by computing in the tubular neighborhood,
\[
\partial_\nu w = \partial_{\nu_i} v_i (\Phi_i) + o(\nabla \eta_i)
\]
To sum up,
\[
\Delta_{\mathbb S^{d-1}} w_i  = -(d-1)w_i + R_i \quad \text{on $\Gamma$}
\]
\[
w_i = 0 \quad\text{on $\partial \Gamma$}
\]
\[
\partial_\nu w_i = 1 + o(\nabla \eta_i) \quad\text{on $\partial \Gamma$}
\]
By the Schauder estimate, $w_i$ is uniformly $C^{2,\alpha}$ bounded on $\overline\Gamma$. Thus, passing to a subsequence, $w_i \rightarrow w_\infty$ in $C^{2}(\overline\Gamma)$ for some $w_\infty \in C^2(\overline\Gamma)$. Then
\[
\Delta_{\mathbb S^{d-1}} w_\infty  = -(d-1)w_\infty \quad \text{on $\Gamma$}
\]
\[
w_\infty = 0 \quad\text{on $\partial \Gamma$}
\]
\[
\partial_\nu w_\infty = 1 \quad\text{on $\partial \Gamma$}
\]
In particular, $w_\infty$, $u := U\mid_{\overline{\Gamma}}$ are eigenfunctions of the first Dirichlet eigenvalue. Therefore, $w_\infty = C u$ for some $C >0$, where the Neumann boundary condition forces $C = 1$. Hence, $w_i \rightarrow u$ in $C^{2}$.\\
Therefore,
\[
\int_{A_{1,2}} |U - V_i|^2 = \int_{1}^2 \int_{\mathbb S^{d-1}} |u(\theta) - v_i(\theta)|^2 r^{d+1}\, d\theta \, dr
\]
converges to 0 as $\int_{\mathbb S^{d-1}} |u(\theta) - w_i(\theta)|^2 \, d\theta$ and $\int_{\mathbb S^{d-1}} |w_i(\theta) - v_i(\theta)|^2 \, d\theta$ converges to 0.
\end{proof}

\section{Non-degenerate weak solutions near a minimizing cone}
As the first step, we consider graphicality results away from the singularity of non-degenerate weak solutions close to one-homogeneous global minimizer. Throughout this section, we assume that $U :\mathbb R^d \rightarrow [0, \infty)$ is a minimizing cone with $\mathrm{Sing}(U) = \{0\}$ and $C_{r_1, r_2}$ is the corresponding tubular neighborhood as defined in Definition \ref{Definition_tubular_neighborhood} with suitably chosen $\tau > 0$ so that $C_{0, \infty}$ is a tubular neighborhood of $\partial \{U > 0\} \setminus \{0\}$.

\begin{lemma} \label{Lemma_graphicality}
Let $\epsilon > 0, k \in \mathbb N$. Then there exists $\delta_0 = \delta_0(d, U, c_1, c_2, \epsilon, k) > 0$ such that if $u$ is a $(c_1, c_2)$-non-degenerate weak solution in $A_{\frac{1}{4}, \frac{3}{4}}$ such that $\|u - U\|_{L^2(A_{\frac{1}{4}, \frac{3}{4}})} < \delta_0$, then
\begin{enumerate}
    \item $u \in C^\omega(\overline{\{u > 0\}}\cap A_{\frac{5}{16}, \frac{1}{2}})$ and $\partial \{u > 0\} \cap A_{\frac{5}{16}, \frac{1}{2}}$ is an analytic hypersurface in $\mathbb R^d$.
    \item there exists $\xi \in C^\omega(\partial \{U > 0\} \cap A_{\frac{5}{16}, \frac{1}{2}})$ such that
    \begin{equation}
    \begin{cases}
        \partial \{ u > 0\} \cap C_{\frac{5}{16}, \frac{1}{2}} = G_{\partial \{U > 0\} \cap A_{\frac{5}{16}, \frac{1}{2}}}(\xi) \cap C_{\frac{5}{16}, \frac{1}{2}} \\
        
        \|\xi\|_{C^k(\partial \{U > 0\} \cap A_{\frac{5}{16}, \frac{1}{2}})} \leq \epsilon \\

        \{u > 0\} \cap C_{\frac{5}{16}, \frac{1}{2}} = \{x + t\nabla U (x) : x \in \partial \{ U > 0\} \cap A_{\frac{5}{16}, \frac{1}{2}}, t > \xi(x) \} \cap C_{\frac{5}{16}, \frac{1}{2}}
    \end{cases}
    \end{equation}
\end{enumerate}
\end{lemma}
\begin{proof}
\cite[Lemma 3.5]{ERZ25}
\end{proof}

The next lemma says that we have graphicality on $A_{\rho, 1}$ when $U$ is integrable through rotations and the difference of Weiss energy at scales $\rho$ and 1 is small.

\begin{lemma} \label{Lemma_extension}
We further assume that $U$ is integrable through rotations. For $\epsilon > 0$, there exists $\delta(U, \epsilon, c_1, c_2) > 0$ such that the following holds:\\
If $u : B_1 \rightarrow [0, \infty)$ is a $(c_1, c_2)$ non-degenerate weak solution and for some $\rho \in [0, \frac{1}{4}]$
\[
    \|u - U\|_{L^2(\partial B_1)} \leq \delta, \quad W_u(0, 1) - W_u(0, \rho/16) \leq \delta
\]
then there exists $\xi \in C^2(\partial\{U > 0\}\cap A_{\rho, \frac{3}{4}})$ such that
\begin{equation}
\begin{cases}
    \partial \{u > 0\} \cap C_{\rho, \frac{3}{4}} = G_{\partial\{U > 0\}\cap A_{\rho, \frac{3}{4}}}(\xi) \cap C_{\rho, \frac{3}{4}} \\

    |x|^{-1}|\xi| + |\nabla \xi| + |x| |\nabla^2 \xi| \leq \epsilon \text{ for $x \in \partial \{ U > 0\}\cap A_{\rho, \frac{3}{4}}$} \\

    \{ u > 0 \} \cap C_{\rho, \frac{3}{4}} = \{x + t\nabla U (x) : x \in \partial \{ U > 0\} \cap A_{\rho, \frac{3}{4}}, t > \xi(x) \} \cap C_{\rho, \frac{3}{4}}
\end{cases}
\end{equation}
\end{lemma}
\begin{proof}
Follows from Lemma \ref{Lemma_graphicality} and \cite[Theorem 4.3]{ERZ25}.
\end{proof}

The next lemma shows that small $L^2$ difference of a cone after a rotation $q$ implies that we may replace $q$ by a rotation $q'$ which is closed to the identity.

\begin{lemma} \label{Lemma_small_rotation}
Let $\epsilon > 0$. There exists $\delta(U, \epsilon) > 0$ such that for $q \in SO(d)$, if $\|U - U(q\cdot)\|_{L^2(A_{1,2})} \leq \delta$, then there exists $q' \in SO(d)$ such that $|q'-\mathrm{Id}| \leq \epsilon$ and $U(q\cdot) \equiv U(q'\cdot)$.
\end{lemma}
\begin{proof}
Suppose not. Then we have $q_i \in SO(d)$ with $\|U - U(q_i\cdot)\|_{L^2(A_{1, 2})} \rightarrow 0$ and for every $q' \in SO(d)$, if $|q'-\mathrm{Id}| \leq \epsilon$, then $U(q\cdot) \not\equiv U(q'\cdot)$. By compactness, passing to a subsequence, $q_i \rightarrow q_\infty$ for some $q_\infty \in SO(d)$. By taking $i \rightarrow \infty$, $U \equiv U(q_\infty\cdot)$. In particular, $U(q_\infty^{-1}q_i\cdot) \equiv U(q_i \cdot)$ and by setting $q'_i := q_\infty^{-1}q_i$, we reach a contradiction.
\end{proof}

Now, we are ready to prove the Theorem \ref{main_theorem}.

\begin{proof}[Proof of Theorem \ref{main_theorem}]
Suppose not. Then we have a sequence $(u_i)$ of $(c_1,c_2, c_3)$-non-degenerate weak solutions so that $u_i \rightarrow U$ in $L^2(B_2)$ and the conclusion fails for $\epsilon > 0$. For each $s \in (0, 1)$, we claim that there are only finitely many $i$ so that $\partial \{ u_i > 0\} \cap B_s = \emptyset$. If not, taking a subsequence, $\partial \{ u_i > 0\} \cap B_s = \emptyset$ for all $i$. In particular, $u_i$ is harmonic in $B_s$. From $L^2$ convergence, $U$ must be harmonic on $B_s$, which is not true.

By Lemma \ref{Lemma_compactness}, \ref{Lemma_positivity_set_convergence}, passing to a subsequence, $u_i \rightarrow U$ in $C(B_2)\cap H^1(B_2)$, $\partial \{ u_i > 0\} \rightarrow \partial \{ u_0 > 0\}$ locally in Hausdorff in $B_2$, and $\chi_{\{u_i > 0\}} \rightarrow \chi_{\{U > 0\}}$ in $L^1(B_2)$. In particular, we have $W_{u_i}(x, r) \rightarrow W_U(x, r)$ for $x \in B_2$ and $r \in (0, d(x, \partial B_2))$.

Fix $\beta_0, \delta_0 > 0$ which will be determined later. For each $i$ and $x \in \overline{B}_{3/2}$, define
\begin{equation} \label{Equation_definition_of_r_i}
\begin{split}
    r_i(x) &:= \inf \{r > 0  : W_{u_i}(x,s) \geq W_U(0) - \beta_0 \text{ for all $s \in (r, 1)$}\}\\
    r_i &:= \inf_{x \in \overline{B}_{3/2}} r_i(x) \text{ with $r_i = r_i(x_i)$}
\end{split}
\end{equation}
where $x_i \in \overline{B}_{3/2}$ exists since $r_i$ is lower semicontinuous.

In particular,
\begin{equation} \label{Equation_Weiss_bound}
    W_{u_i}(x, r_i) \leq W_U(0) - \beta_0 \text{ for any $x \in \overline{B}_{3/2}$}
\end{equation}
when $r_i > 0$.

Passing to a subsequence, we may assume that $x_i \rightarrow x_\infty$ for some $x_\infty \in \overline{B}_{3/2}$. Since $W_{u_i}(0, r) \rightarrow W_U(0)$ for all $r \in (0, 2)$, we have $r_i \rightarrow 0$. We show that
\[
W_{u_i}(x_i, s) \rightarrow W_U(x_\infty, s)
\]
for each $s \in (0, 1/2)$.

It suffices to show that
\[
\int_{\partial B_s(x_i)} u_i^2\, d\mathcal{H}^{d-1} \rightarrow \int_{\partial B_s(x_\infty)} U^2\, d\mathcal{H}^{d-1}
\]
We have
\begin{equation}
\begin{split}
\left| \int_{\partial B_s(x_i)} u_i^2\, d\mathcal{H}^{d-1} - \int_{\partial B_s(x_\infty)} U^2\, d\mathcal{H}^{d-1} \right|
&\leq \left| \int_{\partial B_s(x_i)} u_i^2\, d\mathcal{H}^{d-1} - \int_{\partial B_s(x_\infty)} u_i^2\, d\mathcal{H}^{d-1} \right| \\
&+ \left| \int_{\partial B_s(x_\infty)} u_i^2 - U^2\, d\mathcal{H}^{d-1} \right|    
\end{split}
\end{equation}
The second term follows from the uniform convergence. Consider the first term. By Lemma \ref{Lemma_uniform_lipschitz}, $(u_i^2)$ is uniformly Lipschitz on $B_2$. Hence.
\begin{equation}
\begin{split}
\left| \int_{\partial B_s(x_i)} u_i^2\, d\mathcal{H}^{d-1} - \int_{\partial B_s(x_\infty)} u_i^2\, d\mathcal{H}^{d-1} \right|
&= \left| \int_{\partial B_s(x_\infty)} u_i^2(x - (x_\infty - x_i)) - u_i^2(x)\, d\mathcal{H}^{d-1}(x) \right| \\
&\leq C \mathcal{H}^{d-1}(\partial B_s(x_\infty)) |x_\infty - x_i| \rightarrow 0
\end{split}
\end{equation}

Next, we show that $x_\infty = 0$. For each $s > 0$, we have
\[
    W_U(0) - \beta_0 \leq W_{u_i}(x_i, r_i) \leq W_{u_i}(x_i, s)
\]
for $i \gg 1$ depending on $s$.

Taking limit as $i \rightarrow \infty$ and then $s \downarrow 0$,
\[
    W_U(0) - \beta_0 \leq W_U(x_\infty, 0)
\]
By choosing $\beta_0$ small depending only on $U$ and $\mathbb R^d$, we may assume $W_U(0) - \beta_0 > \frac{\omega_d}{2}$. By Lemma \ref{Lemma_characterization_of_Weiss_for_minimizer}, we conclude that $x_\infty = 0$.

Since $W_{u_i}(x_i, 1) \rightarrow W_U(0)$, for $i \gg 1$ and $W_{(u_i)_{x_i}}(0, 1) = W_{u_i}(x_i, 1)$, we have $W_{(u_i)_{x_i}}(0, 1) \leq W_U(0) + \beta_0$. In particular, $W_{(u_i)_{x_i}}(0, 1) - W_{(u_i)_{x_i}}(0, \rho_i/16) < 2\beta_0$ where 
\begin{equation} \label{Equation_definition_of_rho_i}
    \rho_i := 16r_i.
\end{equation}
Therefore, by Lemma \ref{Lemma_extension} and ensuring that $2\beta_0$ is small depending only on $U, \delta_0$, and $i \gg 1$ so that $\rho_i \leq \frac{1}{4}$ and $\|(u_i)_{x_i} - U\|_{L^2(\partial B_1)}$ small, there exists $\xi_i \in C^2(\partial\{U > 0\}\cap A_{\rho_i, \frac{3}{4}})$ such that
\begin{equation} \label{Equation_u_i_x_i_over_cone}
\begin{cases}
    \partial \{(u_i)_{x_i} > 0\} \cap C_{\rho_i, \frac{3}{4}} = G_{\partial\{U > 0\}\cap A_{\rho_i, \frac{3}{4}}}(\xi_i) \cap C_{\rho_i, \frac{3}{4}} \\
    
    |x|^{-1}|\xi_i| + |\nabla \xi_i| + |x| |\nabla^2 \xi_i| \leq \delta_0 \text{ for $x \in \partial \{ U > 0\} \cap A_{\rho_i, \frac{3}{4}}$} \\
    
    \{(u_i)_{x_i} > 0\} \cap C_{\rho_i, \frac{3}{4}} = \{ x + t \nabla U (x) : x \in \partial \{ U > 0\} \cap A_{\rho_i, \frac{3}{4}}, t > \xi_i(x) \} \cap C_{\rho_i, \frac{3}{4}}.
\end{cases}
\end{equation}

Suppose there are infinitely many $i$ such that $\rho_i = 0$. Passing to a subsequence, assume that $\rho_i = 0$ for all $i$. Then by setting $\bar{\xi}_i(x) := \xi_i(x-x_i)$ for $x \in F_{x_i, 0, \mathrm{Id}} \cap B_{\frac{1}{2}}(x_i)$, where $F_{x_i, 0, \mathrm{Id}}$ is defined in (\ref{Definition_F}),
\begin{equation}
\begin{cases}
    \partial\{ u_i > 0\} \cap B_{\frac{1}{2}}(x_i) = G_{F_{x_i, 0, \mathrm{Id}} \cap B_{\frac{1}{2}}(x_i)}(\bar{\xi}_i) \cap B_{\frac{1}{2}}(x_i)\\

    |x-x_i|^{-1} |\bar{\xi}| + |\nabla \bar{\xi}| + |x-x_i||\nabla^2 \bar{\xi}| \leq \epsilon_0 \text{ for $x \in F_{x_i, 0, \mathrm{Id}} \cap B_{\frac{1}{2}}(x_i)$}
\end{cases}
\end{equation}
Therefore, for $\epsilon_0$ small and $i \gg 1$, we get a contradiction.

Now, suppose $\rho_i > 0$ for all $i$.
Let $v_i := (u_i)_{x_i, \rho_i}$ defined on $B_{\frac{1}{\rho_i}}$, $\xi_i' := (\xi_i)_{\rho_i}$ defined on $\partial\{ U > 0\} \cap A_{1, \frac{3}{4\rho_i}}$. Then
\begin{equation}
\begin{cases}
    \partial \{v_i > 0\} \cap C_{1, \frac{3}{4\rho_i}} = G_{\partial\{ U > 0\} \cap A_{1, \frac{3}{4\rho_i}}}(\xi'_i) \cap C_{1, \frac{3}{4\rho_i}} \\
    
    |x|^{-1} |\xi_i'| + |\nabla\xi_i'| + |x| |\nabla^2\xi_i'| \leq \delta_0 \text{ for $x \in \partial \{U > 0\} \cap A_{1, \frac{3}{4\rho_i}}$} \\
    
    \{v_i > 0\} \cap C_{1, \frac{3}{4\rho_i}} = \{ x + t \nabla U (x) : x \in \partial \{ U > 0\} \cap A_{1, \frac{3}{4\rho_i}}, t > \xi_i'(x) \} \cap C_{1, \frac{3}{4\rho_i}}.
\end{cases}
\end{equation}

For all $i$, $v_i$ is a $(c_1, c_2, c_3)$-non-degenerate weak solution in $B_{\frac{1}{\rho_i}}$ and with non-empty $\partial \{v_i > 0\}$. Hence, by Lemma \ref{Lemma_compactness}, \ref{Lemma_positivity_set_convergence}, upto a subsequence, we have $v_i \rightarrow v_\infty$ locally uniformly in $\mathbb R^d$, $\partial \{ v_i> 0\} \rightarrow \partial\{v_\infty > 0\}$ locally in Hausdorff distance in $\mathbb R^d$, and $\chi_{\{v_i > 0\}} \rightarrow \chi_{\{v_\infty > 0\}}$ in $L^1_{\mathrm{loc}}(\mathbb R^d)$, where $v_\infty : \mathbb R^d \rightarrow \mathbb R$ is a $(c_1, c_2)$-non-degenerate weak solution in $\mathbb R^d$. Moreover, $\partial \{v_\infty > 0\}$ is non-empty.

By Arzela-Asocoli and Lemma \ref{Lemma_extension} (with a minor modification, one can obtain a uniform $C^3$ bound on $(\xi_i')$), upto a subsequence, $\xi_i'\rightarrow \xi_\infty'$ in $C^2_{\mathrm{loc}}(\partial\{U > 0\} \cap A_{1, \infty})$. In particular, $G_{\partial \{U > 0\} \cap A_{1, \frac{1}{2\rho_i}}}(\xi_i') \rightarrow G_{\partial \{U > 0\} \cap A_{1, \infty}}(\xi_\infty')$ locally in Hausdorff distance in $A_{1, \infty}$. Hence, by the uniqueness of Hausdorff limit,
\begin{equation} \label{Equation_v_infty_graph}
\begin{cases}
    \partial \{ v_\infty > 0\} \cap A_{1, \infty} = G_{\partial \{U > 0\} \cap A_{1, \infty}}(\xi_\infty') \cap A_{1, \infty} \\

    |x|^{-1} |\xi_\infty'| + |\nabla \xi_\infty'| + |x| |\nabla^2 \xi_\infty'| \leq \delta_0 \text{ for $x \in \partial \{ U > 0 \} \cap A_{1, \infty}$}
\end{cases}
\end{equation}
Also from the local Hausdorff convergence $\partial \{ v_i> 0\} \rightarrow \partial\{v_\infty > 0\}$ and the $L^1_{\mathrm{loc}}$ convergence $\chi_{\{v_i > 0\}} \rightarrow \chi_{\{v_\infty > 0\}}$, we conclude that
\[
    \{ v_\infty > 0 \} \cap C_{1, \infty} = \{ x + t \nabla U(x) : x \in \partial \{ U > 0 \} \cap A_{1, \infty}, t > \xi_\infty'(x) \} \cap C_{1, \infty}.
\]

By Lemma \ref{Lemma_compactness}, \ref{Lemma_positivity_set_convergence}, there exist $R_j \uparrow \infty$ and a $(c_1, c_2)$-non-degenerate weak solution $V : \mathbb R^d \rightarrow \mathbb R$ such that $(v_\infty)_{R_j} \rightarrow V$ locally uniformly in $\mathbb R^d$. With the same argument as above, there exists $\tilde{\xi} \in C^2(\partial \{ U > 0\}\cap A_{0, \infty})$ such that
\begin{equation}
\begin{cases}
    \partial \{ V > 0\} \cap C_{0, \infty} = G_{\partial \{U > 0\} \cap A_{0, \infty}}(\tilde{\xi}) \cap A_{0, \infty} \\

    |x|^{-1} |\tilde{\xi}| + |\nabla \tilde{\xi}| + |x| |\nabla^2 \tilde{\xi}| \leq \delta_0 \text{ for $x \in \partial \{ U > 0\} \cap A_{0, \infty}$}\\

    \{ V > 0 \} \cap C_{0, \infty} = \{ x + t\nabla U(x) : x \in \partial \{U > 0\} \cap A_{0, \infty}, t > \tilde{\xi}(x)\} \cap C_{0, \infty}
\end{cases}
\end{equation}

Let $\eta, \epsilon' > 0$. By Theorem \ref{Theorem_SMALL_FREE_BOUNDARY_GRAPH_IMPLIES_SMALL_L^2}, choosing $\delta_0(U, \eta)$ small, we have $\|U - V\|_{L^2} \leq \eta$. By \cite[Corollary 4.7]{ERZ25} and choosing $\eta(U)$ small, there exists $q \in SO(d)$ such that $V = U(q^{-1}\cdot)$. By Lemma \ref{Lemma_small_rotation}, choosing $\eta(U, \epsilon')$ small, we may pick $q \in SO(d)$ such that $|q- \mathrm{Id}| \leq \epsilon'$. By the rigidity Theorem \ref{Theorem_blowdown_rigidity}, there exists $a \in \mathbb R^d$ and $t > 0$ such that $v_\infty$ coincides with $U(q^{-1}(\cdot - a)), \overline{U}_{\frac{1}{t}}(q^{-1}(\cdot - a))$ or $\underline{U}_{\frac{1}{t}}(q^{-1}(\cdot - a))$.

By (\ref{Equation_definition_of_r_i}), (\ref{Equation_Weiss_bound}), (\ref{Equation_definition_of_rho_i}) and taking limit as $i \rightarrow \infty$, we have
\begin{equation} \label{Equation_Weiss_bound_on_v_infty}
    W_{v_\infty}\left(x, 1/16\right) \leq W_U(0) - \beta_0 \text{ for all $x \in \overline{B}_{3/2}$.}
\end{equation}

We show that $v_\infty$ cannot be $U(q^{-1}(\cdot - a))$. Suppose not. Since $a$ is a singular point of $\partial \{v_\infty > 0\}$ and $\{v_\infty > 0\} \cap C_{1, \infty}$ is regular, provided that $\delta_0(U)$ is small, we have $a \in \overline{B}_{3/2}$. Also, notice that
\[
    W_U(0) = W_{U(q^{-1}\cdot)}(0, 1/16) = W_{v_\infty}(a, 1/16)
\]
Hence, from (\ref{Equation_Weiss_bound_on_v_infty}) with $x = a$, we get a contradiction.

Let $\epsilon''>0$. Since the unique tangent cone of $v_\infty$ is $U(q^{-1}\cdot)$, there exists $R_0 > 0$ and $f \in C^2(F_{0, q, 0} \cap A_{\frac{R_0}{2}, \infty})$ such that
\begin{equation}
\begin{cases}
    F_{a, q, t} \cap A_{R_0, \infty} = G_{F_{0, q, 0} \cap A_{\frac{R_0}{2}, \infty}}(f) \cap A_{R_0, \infty} \\

    |x|^{-1}|f| + |\nabla f| + |x| |\nabla^2 f| \leq \epsilon'' \text{ for $x \in F_{0, q, 0} \cap A_{\frac{R_0}{2}, \infty}$}.
\end{cases}
\end{equation}
Let $f_i := f_{\frac{1}{\rho_i}}$. Then
\begin{equation} \label{Equation_F_over_rotated_cone}
\begin{cases}
    F_{\rho_i a, q, \rho_i t} \cap A_{\rho_i R_0, \infty} = G_{F_{0, q, 0} \cap A_{\frac{\rho_iR_0}{2}, \infty}}(f_i) \cap A_{\rho_i R_0, \infty} \\

    |x|^{-1}|f_i| + |\nabla f_i| + |x| |\nabla^2 f_i| \leq \epsilon'' \text{ for $x \in F_{0, q, 0} \cap A_{\frac{\rho_i R_0}{2}, \infty}$}.
\end{cases}
\end{equation}

Also, by choosing $\epsilon'(U, \epsilon'')$ small, there exists $g \in C^2(\partial \{ U > 0\} \cap A_{\frac{R_0}{2}, \infty})$ such that
\begin{equation}
\begin{cases}
    F_{0, q, 0} \cap A_{R_0, \infty} = G_{\partial \{ U > 0\} \cap A_{\frac{R_0}{2}, \infty}}(g) \cap A_{R_0, \infty} \\

    |x|^{-1}|g| + |\nabla g| + |x| |\nabla^2 g| \leq \epsilon'' \text{ for $x \in \partial \{ U > 0\} \cap A_{\frac{R_0}{2}, \infty}$}.
\end{cases}
\end{equation}

Let $g_i := g_{\frac{1}{\rho_i}}$. Then
\begin{equation} \label{Equation_rotated_cone_over_cone}
\begin{cases}
    F_{0, q, 0} \cap A_{\rho_i R_0, \infty} = G_{\partial \{ U > 0\} \cap A_{\frac{\rho_i R_0}{2}, \infty}}(g_i) \cap A_{\rho_i R_0, \infty} \\

    |x|^{-1}|g_i| + |\nabla g_i| + |x| |\nabla^2 g_i| \leq \epsilon'' \text{ for $x \in \partial \{ U > 0\} \cap A_{\frac{\rho_i R_0}{2}, \infty}$}.
\end{cases}
\end{equation}

From (\ref{Equation_F_over_rotated_cone}), (\ref{Equation_rotated_cone_over_cone}) and choosing $\epsilon'' (U)$ small, there exists $C(U) > 0$ and $h_i \in C^2(\partial \{ U > 0\} \cap A_{\rho_i R_0, 2})$ such that
\begin{equation} \label{Equation_F_over_cone}
\begin{cases}
    F_{\rho_i a, q, \rho_i t} \cap A_{2\rho_i R_0, 1} = G_{\partial \{ U > 0\} \cap A_{\rho_i R_0, 2}}(h_i) \cap A_{2\rho_i R_0, 1} \\

    |x|^{-1}|h_i| + |\nabla h_i| + |x| |\nabla^2 h_i| \leq C\epsilon'' \text{ for $x \in \partial \{ U > 0\} \cap A_{\rho_i R_0, 2}$}.
\end{cases}
\end{equation}

From (\ref{Equation_u_i_x_i_over_cone}), (\ref{Equation_F_over_cone}) and choosing $\epsilon''(U), \delta_0(U)$ small such that $\delta_0 \leq \epsilon''$, there exists $p_i \in C^2(F_{\rho_i a, q, \rho_i t} \cap A_{4\rho_i R_0, \frac{5}{8}})$ such that
\begin{equation}
\begin{cases}
    \partial \{ (u_i)_{x_i} > 0 \} \cap A_{5\rho_i R_0, \frac{1}{3}} = G_{F_{\rho_i a, q, \rho_i t} \cap A_{4\rho_i R_0, \frac{5}{8}}}(p_i) \cap A_{5\rho_i R_0, \frac{1}{3}} \\

    |x|^{-1}|p_i| + |\nabla p_i| + |x| |\nabla^2 p_i| \leq C\epsilon'' \text{ for $x \in F_{\rho_i a, q, \rho_i t} \cap A_{4\rho_i R_0, \frac{5}{8}}$}.
\end{cases}
\end{equation}

For $x \in A_{4\rho_i R_0, \frac{5}{8}}$, we have $\rho_i \leq \frac{|x|}{4R_0}$ and so
\[
    \left(1 - \frac{|a|}{4R_0} \right) |x| \leq |x - \rho_i a| \leq \left(1 + \frac{|a|}{4R_0} \right) |x|
\]
where we also impose that $R_0 \geq \max \{1, |a|\}$.

Therefore, by increasing $C$,
\begin{equation} \label{Equation_u_i_x_i_over_F_outside}
\begin{cases}
    \partial \{ (u_i)_{x_i} > 0 \} \cap A_{5\rho_i R_0, \frac{1}{3}} = G_{F_{\rho_i a, q, \rho_i t} \cap A_{4\rho_i R_0, \frac{5}{8}}}(p_i) \cap A_{5\rho_i R_0, \frac{1}{3}} \\

    |x-\rho_i a|^{-1}|p_i| + |\nabla p_i| + |x - \rho_i a| |\nabla^2 p_i| \leq C\epsilon'' \text{ for $x \in F_{\rho_i a, q, \rho_i t} \cap A_{4\rho_i R_0, \frac{5}{8}}$}.
\end{cases}
\end{equation}

Now, wlog suppose that $v_\infty = \overline{U}_{\frac{1}{t}}(q^{-1}(\cdot - a))$. For $i \gg 1$, there exists $\eta_i \in C^2(F_{a, q, t} \cap B_{7R_0})$ such that
\begin{equation}
    \partial \{  v_i > 0 \} \cap B_{6R_0} = G_{F_{a, q, t} \cap B_{7R_0}}(\eta_i) \cap B_{6R_0}
\end{equation}
\begin{equation}
    \|\eta_i\|_{C^2(F_{a, q, t} \cap B_{7R_0})} \rightarrow 0
\end{equation}

For $x \in F_{a, q, t} \cap B_{7R_0}$, we have $|x-a| \geq t$ and so
\begin{equation}
    |x-a|^{-1} |\eta_i| + |\nabla \eta_i| + |x-a| | \nabla^2 \eta_i| \leq C(t, R_0) \|\eta_i\|_{C^2(F_{a, q, t} \cap B_{7R_0})}
\end{equation}

In particular, for $i \gg 1$, we have
\begin{equation}
    |x-a|^{-1} |\eta_i| + |\nabla \eta_i| + |x-a| | \nabla^2 \eta_i| \leq \epsilon''
\end{equation}

Let $\tilde{\eta}_i := (\eta_i)_{\frac{1}{\rho_i}}$. Then
\begin{equation} \label{Equation_u_i_x_i_over_F_inside}
\begin{cases}
    \partial \{ (u_i)_{x_i} > 0 \} \cap B_{6\rho_i R_0} = G_{F_{\rho_i a, q, \rho_i t} \cap B_{7\rho_i R_0}}(\tilde{\eta}_i) \cap B_{6\rho_i R_0} \\

    |x - \rho_i a|^{-1} |\tilde{\eta}_i| + |\nabla \tilde{\eta}_i| + |x - \rho_i a||\nabla^2 \tilde{\eta}_i| \leq \epsilon'' \text{ for $x \in F_{\rho_i a, q, \rho_i t} \cap B_{7\rho_i R_0}$}.
\end{cases}
\end{equation}

Let $a_i := x_i + \rho_i a$ and $t_i := \rho_i t$. From (\ref{Equation_u_i_x_i_over_F_outside}), (\ref{Equation_u_i_x_i_over_F_inside}), we have $\zeta_i \in C^2(F_{a_i, q, t_i} \cap B_{\frac{5}{8}}(x_i))$ such that
\begin{equation}
\begin{cases}
    \partial \{ u_i > 0 \} \cap B_{\frac{1}{3}}(x_i) = G_{F_{a_i, q, t_i} \cap B_{\frac{5}{8}}(x_i)}(\zeta_i) \cap B_{\frac{1}{3}}(x_i) \\

    |x - a_i|^{-1}|\zeta_i| + |\nabla \zeta_i| + |x-a_i| |\nabla^2 \zeta_i| \leq C \epsilon'' \text{ for $x \in F_{a_i, q, t_i} \cap B_{\frac{5}{8}}(x_i)$}.
\end{cases}
\end{equation}
Therefore, ensuring $C\epsilon'' \leq \min\{\epsilon, 1/2\}$, $\epsilon' \leq \epsilon/2$, for $i \gg 1$, since $B_{\frac{1}{2}}(a_i) \subseteq B_{\frac{5}{8}}(x_i)$ and $B_{\frac{1}{4}}(a_i) \subseteq B_{\frac{1}{3}}(x_i)$, we get a contradiction.
\end{proof}

\bibliographystyle{amsalpha}
\bibliography{references}

\end{document}